\documentclass[11pt, reqno, a4paper]{amsart}

\usepackage{amsfonts}
\usepackage[T1]{fontenc}
\usepackage{enumerate}
\usepackage{amsmath}
\usepackage[latin1]{inputenc}
\usepackage{amssymb}
\usepackage{verbatim}
\usepackage{amsthm}
\usepackage{hyperref}
\usepackage[arrow, matrix, curve]{xy}
\usepackage{tikz-cd}
\usepackage{tikz}
\usetikzlibrary{arrows,matrix}
\usepackage{relsize}
\usepackage{MnSymbol}

\usepackage{listings}
\usepackage{color}
\definecolor{Grau}{gray}{.94}

\usepackage{geometry}
\DeclareMathOperator{\fpt}{fpt}
\DeclareMathOperator{\lct}{lct}
\DeclareMathOperator{\Proj}{Proj}
\DeclareMathOperator{\ord}{ord}
\DeclareMathOperator{\ggT}{gcd}
\DeclareMathOperator{\h}{ht}
\DeclareMathOperator{\codim}{codim}

\usepackage{dsfont}

\usepackage{mathrsfs}
\usepackage{stmaryrd}

\address{Susanne Müller\\ Johannes Gutenberg-Universität Mainz \\ Fachbereich 08 \\ Staudingerweg 9 \\ 55099 Mainz \\ Germany}
\email{susanne.mueller@uni-mainz.de}

\begin{document} 

\title{$F$-pure threshold and height of quasi-homogeneous polynomials}
\author{Susanne Müller}
\begin{abstract}
\small{We consider a quasi-homogeneous polynomial $f \in \mathbb{Z}[x_0, \ldots, x_N]$ of degree $w$ equal to the degree of $x_0 \cdots x_N$ and show that the $F$-pure threshold of the reduction $f_p \in \mathbb{F}_p[x_0, \ldots, x_N]$ is equal to the log canonical threshold if and only if the height of the Artin-Mazur formal group associated to $H^{N-1}\left( X, {\mathbb{G}}_{m,X} \right)$, where $X$ is the hypersurface given by $f$, is equal to 1.
We also prove that a similar result holds for Fermat hypersurfaces of degree $>N+1$.
Furthermore, we give examples of weighted Delsarte surfaces which show that other values of the $F$-pure threshold of a quasi-homogeneous polynomial of degree $w$ cannot be characterized by the height.}
\end{abstract}
\maketitle

\newtheoremstyle{Remark}
  {5pt}                   
  {5pt}                   
  {\normalfont}           
  {}                      
  {\bfseries}             
  {.}                     
  {2mm}                   
  {}                      

\newtheorem{definition}{Definition}[section]
\newtheorem{theorem}[definition]{Theorem}
\newtheorem{proposition}[definition]{Proposition}
\newtheorem{lemma}[definition]{Lemma}
\newtheorem{corollary}[definition]{Corollary}

\newtheorem*{theorem_o}{Theorem} 
\newtheorem*{lemma_o}{Lemma} 

\theoremstyle{Remark}
\newtheorem{remark}[definition]{Remark}
\newtheorem{example}[definition]{Example}

\newtheorem*{example_o}{Example}

\section{Introduction}

To any polynomial $f \in \mathbb{F}_p[x_0, \ldots, x_N]$ one can attach an invariant called the $F$-pure threshold, first defined in \cite{TakWat}, \cite{MTW}.  
The $F$-pure threshold, which is a rational number (see \cite{BMS}), is a quantitative measure of the severity of the singularity of $f$. Smaller values of the $F$-pure threshold correspond to a "worse" singularity.
For a short introduction to the theory of $F$-pure thresholds see \cite{MTW} or \cite{Mueller}.
In \cite{Mueller} we proved that for a quasi-homogeneous polynomial $f \in \mathbb{F}_p[x_0, \ldots, x_N]$ of degree $w= \alpha_0+ \ldots + \alpha_N$, where $\alpha_i=\deg(x_i)$, with an isolated singularity and with $p \geq w(N-2)+1$ one has $\fpt(f)=1-\frac{a}{p}$. Here, the integer $0 \leq a \leq N-1$ is the order of vanishing of the Hasse invariant on a certain deformation space of $X=\Proj(R/fR) \subset \mathbb{P}^N \left(\alpha_0, \ldots,  \alpha_N \right)$.

The $F$-pure threshold is the characteristic $p$ analogue of the log canonical threshold $\lct$ in characteristic $0$, which is defined via resolution of singularities. In general, it is difficult to compute the log canonical threshold, but for a quasi-homogeneous polynomial of degree $d$ in $N+1$ variables with an isolated singularity, one can show that
$\lct(f)= \frac{w}{d}$ if  $d \geq w$ and $\lct(f)=1$ otherwise (see \cite{Lazarsfeld}). 
Comparing the log canonical threshold of a polynomial $f \in \mathbb{Z}[x_0, \ldots, x_N]$ with the $F$-pure threshold of its reduction $f_p \in \mathbb{F}_p[x_0, \ldots, x_N]$ it turns out that $\fpt(f_p) \leq \lct(f)$ for all $p$ and $\lim\limits_{p \rightarrow \infty} \fpt(f_p) = \lct(f)$ (\cite{TakWat}, \cite{MTW}).
Furthermore, it is conjectured that for infinitely many primes $p$ one has $\fpt(f_p)=\lct(f)$. But this is wide open.

On the other hand, for a polynomial $f \in \mathbb{Z}\left[x_0, \ldots, x_N\right]$ one can consider the hypersurface $X$ in $\mathbb{P}^N_{\mathbb{Z}}\left(\alpha_0, \ldots,  \alpha_N \right)$ given by $f$ and compute the height of the so-called Artin-Mazur formal group associated to $H^{N-1}\left( X, {\mathbb{G}}_{m,X} \right)$, which is either infinite or an integer greater or equal to $1$. This is another important invariant, uniquely characterizing $1$-dimensional formal groups over an algebraically closed field of positive characteristic by Lazard \cite{Lazard}. 

The aim of this paper is to clarify the connection between the $F$-pure threshold and the height by establishing the following two results. Their proofs will occupy section 3.

\begin{theorem_o}[see Theorem \ref{TheoremCYhom}]\label{TheoremCYhom}
Let $\mathbb{Z}[x_0, \ldots, x_N]$ be the graded polynomial ring with $\alpha_i:=\deg(x_i)$ and set $w:= \alpha_0+  \ldots + \alpha_N$.
Let $f \in \mathbb{Z}[x_0, \ldots, x_N]$ be a quasi-homogeneous polynomial of degree $w$ and type $\alpha=\left( \alpha_0, \ldots, \alpha_N \right)$ with an isolated singularity such that the greatest common divisor of all coefficients of $f$ is $1$. Furthermore, let $X$ be the hypersurface in $\mathbb{P}^N_{\mathbb{Z}}\left( \alpha \right)$ defined by $f$. Let $f_p \in \mathbb{F}_p[x_0, \ldots, x_N]$ be the reduction of $f$ modulo $p$ and assume that $p \geq w(N-2)+1$.
Then $\fpt(f_p)=1=\lct(f)$ if and only if $\h\left( H^{N-1}\left( X , {\mathbb{G}}_{m,X}\right) \right)=1$.
\end{theorem_o}

Furthermore, we show that a similar result holds for Fermat hypersurfaces of degree $>N+1$  :

\begin{lemma_o}
Let $f=x_0^d+ \ldots + x_N^d \in \mathbb{Z}[x_0, \ldots, x_N]$ with $d=N+k$ for $k \geq 2$ and such that $N \geq 2(k-1)$. 
Furthermore, let $d \not \equiv 0 \mod p$. \\
Then $H^{N-1} \left( X, {\mathbb{G}}_{m,X} \right)$ is a direct sum of formal groups of dimension $1$, which are all of height $1$ if and only if $\fpt(f_p)=\lct(f)=\frac{N+1}{d}$.
\end{lemma_o}

We will see that the above statements mean that the $F$-pure threshold is equal to the log canonical threshold if and only if the height of the corresponding Artin-Mazur formal group is equal to its dimension. Since $\fpt(f_p) \leq \lct(f)$, this means that the $F$-pure threshold is equal to its greatest possible value if and only if the height is equal to its smallest possible value.
We suspect that this could hold more generally for quasi-homogeneous polynomials. 
All computations of the height and the $F$-pure threshold in concrete examples support this.

The last part of this paper is dedicated to the following:
Theorem \ref{TheoremCYhom} yields that  
for the integer $a$ from above, $a=0$ holds if and only if $\h \left( H^{N-1}\left( X , {\mathbb{G}}_{m,X}\right) \right)=1$.
Therefore, it is natural to ask whether the other possible values of the $F$-pure threshold (i.e. $1 \leq a \leq N-1$) can also be characterized by $\h\left(H^{N-1}\left( X , {\mathbb{G}}_{m,X}\right)\right)$.
However, we will give two examples of weighted Delsarte surfaces which show that the answer to this question is negative. The first example will have the same height but different $F$-pure threshold and the second one will have the same $F$-pure threshold but the height will differ for two different primes $p$.

\textbf{Acknowledgements.} 
I thank Manuel Blickle and Axel Stäbler for useful discussions and a careful reading of earlier versions of this article.
Furthermore, I thank Duco van Straten for the inspiration to work on this subject and Masha Vlasenko for her valuable advice while familiarizing myself with formal groups.
The author was supported by SFB/Transregio 45 Bonn-Essen-Mainz financed by Deutsche Forschungsgemeinschaft.

\section{Formal groups}

As a preparation for the remainder of the paper we begin with a short introduction to formal groups.
In the theory of formal groups one can choose the point of view of formal power series or the point of view of functors -
we will sketch both in what follows.
For further information about the point of view of formal power series we refer the reader to \cite{Froehlich}, \cite{Hazewinkel}, \cite{Honda} and \cite{Vlasenko}. In \cite{Stienstra} and \cite{Zink} the authors also treat the point of view of functors.

Let $x=(x_1, \ldots, x_m)$ and $y=(y_1, \ldots, y_m)$ be two sets of $m$ variables.
An $m$-dimensional \textbf{formal group law} over a commutative ring $R$ with identity element is an $m$-tuple of power series $F(x,y)=\left( F_1(x,y), \ldots, F_m(x,y) \right)$ with $F_i(x,y) \in R \lsem x,y \rsem$, such that
$$F(x,F(y,z))=F(F(x,y),z) \text{ and }$$
$$F(x,y) \equiv x+y \mod \deg \geq 2.$$
A formal group law is called commutative, if one has in addition that $F_i(x,y)=F_i(y,x)$ for all $i$.

Let $F$ and $G$ be two formal group laws over $R$ of dimension $m_F$ and $m_G$ respectively. A \textbf{homomorphism} $F(x,y) \rightarrow G(x,y)$ over $R$ is an $m_G$-tuple of power series $\varphi$ in $m_F$ variables, such that $\varphi(x) \equiv 0 \mod \deg \geq 1$ and $\varphi \left( F(x,y) \right) = G \left( \varphi(x), \varphi(y) \right)$.
The homomorphism $\varphi(x)$ is an isomorphism if there exists a homomorphism $\psi(x): G(x,y) \rightarrow F(x,y)$ such that $\varphi(\psi(x))=x$ and $\psi(\varphi(x))=x$. 
The morphism $\varphi(x)$ is said to be a strict isomorphism if $\varphi(x) \equiv x \mod \deg \geq 2$.

If $R$ is a ring of characteristic zero, 
then every $m$-dimensional commutative formal group law $F(x,y)$ over $R$ determines a unique $m$-tuple $l(\tau)=\left(l_1(\tau), \ldots, l_m(\tau) \right)$ of power series  in an $m$-tuple of variables $\tau=\left( \tau_1, \ldots, \tau_m\right)$ with coefficients in $R \otimes \mathbb{Q}$ such that
$$l(\tau) \equiv \tau \mod \deg \geq 2 \text{ and }$$ 
$$F(x,y)=l^{-1} \left( l(x)+l(y)\right).$$
This $m$-tuple $l(\tau)$ is called the \textbf{logarithm} of the formal group law $F(x,y)$. 
In the $1$-dimensional case one can write 
$$l(\tau)=\tau+\mathlarger{\mathlarger{\sum}}_{m=2}^{\infty} \frac{b_{m-1}}{m} \tau^m$$ 
with $b_{m-1} \in R$. The name of the logarithm comes from the following example:

\begin{example}
We consider the $1$-dimensional additive formal group law $\mathbb{G}_a$ and the $1$-dimensional multiplicative formal group law $\mathbb{G}_m$, which are both defined over $\mathbb{Z}$. The additive formal group law  
is given by $\mathbb{G}_a(x,y)=x+y$ with logarithm $l(\tau)=\tau$. The multiplicative formal group law  
is given by $\mathbb{G}_m(x,y)=x+y+xy$ and the logarithm is $l(\tau)=\log(1+\tau)=\sum_{n \geq 1} (-1)^{n+1} \frac{1}{n} \tau^n$.
\end{example}

Now, let $F(x,y)$ be an $m$-dimensional formal group law over a field $k$ of characteristic $p>0$.
An important invariant of the formal group law is the \textbf{height} $\h=\h(F)$.
Consider the multiplication by $p$ endomorphism, which is given by
$$[p]_{F}(x)=\underbrace{x+_{F} x+_{F}+ \ldots +_{F} x}_{p \text{ times}}$$
and write $[p]_{F}(x)=\left( H_1(x), \ldots, H_m(x) \right)$. We say that $F(x,y)$ is of finite height, if the ring $k \lsem x_1, \ldots, x_m \rsem$ is a finitely generated module over the subring $k \lsem H_1(x), \ldots, H_m(x) \rsem$. In this case, $k \lsem x_1, \ldots, x_m \rsem$ is free of rank $p^r$, $r \in \mathbb{N}$ over $k \lsem H_1(x), \ldots, H_m(x) \rsem$ and $\h(F):=r$ is called the height of $F(x,y)$ (see \cite[18.3.8]{Hazewinkel}).
If $R$ is a local ring of characteristic zero with residue field $k$ of characteristic $p>0$ and $F(x,y)$ is an $m$-dimensional formal group law over $R$, then we define the height of $F(x,y)$ as the height of the reduction $\overline{F}(x,y)$ of $F(x,y)$ over $k$.

If $F(x,y)$ is a one-dimensional formal group law over a field $k$ of characteristic $p>0$, then this definition says the following:
Let $[p]_{F}(x)$ be the multiplication by $p$ as above. Then one can show (see \cite[18.3.1]{Hazewinkel}) that either $[p]_{F}(x)=0$ or there is a power $q=p^r$ of $p$ such that $[p]_{F}(x)=\beta(x^q)$, $\beta(x) \not\equiv 0 \mod \deg \geq 2$. Then $\h(F)=\infty$ iff $[p]_{F}(x)=0$ and $\h(F)=r$ if $q=p^r$ is the highest power of $p$ such that $[p]_{F}(x)=\beta(x^q)$.

\begin{lemma} \label{direkteSummefG}
Let $F= G \times H$ be a formal group, which is the product of two formal groups $G$ and $H$ of finite heights $\h_G$ respectively $\h_H$. Then $F$ has height $\h_G+\h_H$.
\end{lemma}

\begin{proof}
Write $[p]_{\overline{G}}(x)=\left( G_1(x), \ldots, G_m(x) \right)$ and $[p]_{\overline{H}}(y)=\left( H_1(y), \ldots, H_n(y) \right)$. Since $G$ has height $\h_G$, we know that $k \lsem x_1, \ldots, x_m \rsem$ is a finitely generated module over the subring $k \lsem G_1(x), \ldots, G_m(x) \rsem$ of rank $p^{\h_G}$ and since $H$ has height $\h_H$, we know that $k \lsem y_1, \ldots, y_n \rsem$ is a finitely generated module over the subring $k \lsem H_1(y), \ldots, H_n(y) \rsem$ of rank $p^{\h_H}$. Therefore, $k \lsem x_1, \ldots, x_m, y_1, \ldots , y_n \rsem$ is a finitely generated module over the subring $k \lsem G_1(x), \ldots, G_m(x), H_1(y), \ldots, H_n(y) \rsem$ of rank $p^{\h_G} p^{\h_H}=p^{\h_G+\h_H}$.
\end{proof}

The importance of the height becomes clear by the following classification result:

\begin{theorem}[\cite{Lazard}] 
Let $k$ be an algebraically closed field of positive characteristic.
\begin{enumerate}
	\item For every integer $h \geq 1$ and for $h = \infty$ there exists a $1$-dimensional formal group law of height $h$ over $k$.
	\item Two $1$-dimensional formal group laws over $k$ are isomorphic if and only if they have the same height.
\end{enumerate}
\end{theorem}

Now, let us come to the point of view of functors. 
For this, let $\mathfrak{Nilalg}_R$ denote the category of nil-$R$-algebras, i.e. of $R$-algebras in which every element is nilpotent. The formal affine $m$-space over $R$ is defined as the functor 
$$\mathbb{A}^m_R : \mathfrak{Nilalg}_R \rightarrow \mathfrak{Sets},$$
which sends a nil-$R$-algebra $N$ to the set $N^{(m)}:=N \oplus \cdots \oplus N$ with $m$ factors and which sends a morphism $f$ to the map $f \times \cdots \times f$. 
An $m$-dimensional \textbf{formal group} over $R$ is a functor
$$F: \mathfrak{Nilalg}_R \rightarrow \mathfrak{Abelian} \ \mathfrak{Groups},$$
such that $V \circ F \cong \mathbb{A}^m_R$, where $V: \mathfrak{Abelian} \ \mathfrak{Groups} \rightarrow \mathfrak{Sets}$ is the forgetful functor.
One can show that given a commutative formal group law $F(x,y) \in R \lsem x,y \rsem$ one can associate to $F$ a functor $F: \mathfrak{Nilalg}_R \rightarrow \mathfrak{Abelian} \ \mathfrak{Groups}$, where the group structure is given by the power series $F$. Conversely, given a functor $F: \mathfrak{Nilalg}_R \rightarrow \mathfrak{Abelian} \ \mathfrak{Groups}$, then $F$ is defined  by a formal group law (see \cite{Zink}).

Now, if $F$ is a formal group over $R$, $X$ a scheme over $R$ and $i \in \mathbb{N}_0$, then one can construct the following diagram:
\begin{center}
\begin{tikzcd} 
  \mathfrak{Nilalg}_R \arrow{r}{\mathcal{O}_X \otimes_R \underline{\ \ }} \arrow[end anchor=real west]{rd}{{\mathbb{G}}_{m,\mathcal{O}_X}} \arrow[end anchor=real west, swap]{rdd}{H^i \left( X, {\mathbb{G}}_{m,\mathcal{O}_X} \right)} &[4em]
 \
 \mathfrak{Sheaves} \ \mathfrak{of} \ \mathfrak{nil}-R-\mathfrak{algebras} \ \mathfrak{on} \ X \arrow{d}{{\mathbb{G}}_m} \\
			 & \mathfrak{Sheaves} \ \mathfrak{of} \ \mathfrak{abelian} \  \mathfrak{groups} \ \mathfrak{on} \ X \arrow{d}{H^i} \\
       & \mathfrak{Abelian} \ \mathfrak{Groups} 
\end{tikzcd}
\end{center}
Here $\mathcal{O}_X \otimes_R \underline{\ \ }$ assigns to a nil-$R$-algebra $A$ the sheaf $\mathcal{O}_X \otimes_R A$ associated with the pre-sheaf $U \mapsto \Gamma(U,\mathcal{O}_X) \otimes_R A$ for $U$ open. The functor ${\mathbb{G}}_m$ assigns to a sheaf $\mathfrak{a}$ of nil-$R$-algebras on $X$ the sheaf of abelian groups ${\mathbb{G}}_m(\mathfrak{a})$ defined by $\Gamma(U,{\mathbb{G}}_m(\mathfrak{a}))={\mathbb{G}}_m(\Gamma(U,\mathfrak{a}))$ for $U \subset X$ open. The functor $H^i$ is taking $i$-th cohomology and the functors ${\mathbb{G}}_{m, \mathcal{O}_X}$ and $H^i \left( X, {\mathbb{G}}_{m, \mathcal{O}_X} \right)$ are defined by the commutativity of the above diagram.
 Writing ${\mathbb{G}}_{m,X}$ instead of ${\mathbb{G}}_{m, \mathcal{O}_X}$, the functors 
$$H^i \left( X, {\mathbb{G}}_{m,X} \right) : \mathfrak{Nilalg}_R \rightarrow \mathfrak{Abelian} \ \mathfrak{Groups}$$ 
are called \textbf{Artin-Mazur functors}. These functors are not necessarily formal groups but Artin and Mazur (see \cite{ArtinMazur}) give a criterion for $H^{i} \left( X, {\mathbb{G}}_{m,X} \right)$ to be a formal group. The functors $H^1 \left( X, {\mathbb{G}}_{m,X} \right)$ and $H^2 \left( X, {\mathbb{G}}_{m,X} \right)$ are called the formal Picard group and the formal Brauer group, respectively (at least if they are formal groups).

\begin{example}\label{StienstraEx}
In the following, we will often use a criterion of Stienstra (see \cite[Theorem 1]{Stienstra}) for $H^{\bullet} \left( X, {\mathbb{G}}_{m,X} \right)$ to be a formal group:\\
Let $K$ be a noetherian ring and let $X$ be a subscheme of $\mathbb{P}^N_{K}$ defined by the ideal $(F_1, \ldots, F_r)$, where $F_1, \ldots, F_r$ is a regular sequence of homogeneous polynomials in $K[x_0, \ldots, x_N]$. Let $d_i=\deg(F_i)$ and $d = \sum\limits_{i=1}^r d_i$. If $X$ is flat over $K$ and $d_i \geq d-N \geq 1$ for all $i$ then $H^{N-r}\left( X, {\mathbb{G}}_{m,X} \right)$ is a formal group over $K$ of dimension $\binom{d-1}{N}$.

Furthermore, Stienstra computes the logarithm of this formal group. For this, assume that $K$ is flat over $\mathbb{Z}$ and set
$$J:= \left\{i=\left( i_0, \ldots, i_N \right) \in \mathbb{Z}^{N+1} \big| i_0, \ldots, i_N \geq 1, i_0+\ldots+i_N=d\right\}.$$
Then there is a formal group law for $H^{N-r}\left( X, {\mathbb{G}}_{m,X} \right)$ whose logarithm is the tuple $\left( l_i(\tau) \right)_{i \in  J}$ of power series in $\tau=\left( \tau_i \right)_{i \in J}$ with
$$l_i(\tau)= \mathlarger{\mathlarger{\sum}}_{m \geq 1} \mathlarger{\mathlarger{\sum}}_{j \in J} \frac{b_{m-1,i,j}}{m} \tau_j^m,$$
where
$$b_{m-1,i,j}= \text{coefficient of } x_0^{mj_0-i_0} \cdots x_N^{mj_N-i_N} \text{ in } \left(F_1 \cdots F_r \right)^{m-1}.$$
\end{example}

\section{Connection between the $F$-pure threshold and the height}

In order to prove the main theorem of this paper, we first need the following result:

\begin{lemma} \label{VlasenkoCorThm2}
Let $R$ be the ring of integers of a complete absolutely unramified discrete valuation field of characteristic zero and residue characteristic $p>0$, equipped with a lift of the $p$-th power Frobenius on the residue field $R/pR$.
Let $F(x,y) \in R \lsem x,y \rsem$ be a formal group law of dimension $1$ with logarithm
$$l(\tau)=\mathlarger{\mathlarger{\sum}}\limits_{m=1}^{\infty} \frac{b_{m-1}}{m} \tau^m,$$ where $\left\{b_m\right\}_{m \geq 0}$ is a sequence of elements of $R$ with $b_0=1$. Then $\h(F)=1$ if and only if $\ord_p \left( b_{p-1} \right)=0$.
\end{lemma}

\begin{proof}
First, let $\h(F)=1$. Then, by Theorem 2(i) of \cite{Vlasenko}, we get $\ord_p \left( b_{p-1} \right)=1-\left\lfloor \frac{1}{1}\right\rfloor=0$. \\
For the opposite direction, let $\h(F) \neq 1$. Then we have two cases. The first case is $\h(F) = \infty$, which yields $\ord_p \left( b_{p-1} \right) \geq 1$ by Theorem 2(i) of \cite{Vlasenko}. The second case is $\h(F) < \infty$ and $\h(F) \neq 1$. Then, again by Theorem 2(i) of \cite{Vlasenko}, we  conclude that $\ord_p \left( b_{p-1} \right) \geq 1- \left\lfloor \frac{1}{\h(F)}\right\rfloor=1$, since $\h(F)>1$.
\end{proof}

Now, we can prove the main theorem.

\begin{theorem}\label{TheoremCYhom}
Let $\mathbb{Z}[x_0, \ldots, x_N]$ be the graded polynomial ring with $\alpha_i:=\deg(x_i)$ and set $w:= \alpha_0+ \ldots + \alpha_N$.
Let $f \in \mathbb{Z}[x_0, \ldots, x_N]$ be a quasi-homogeneous polynomial of degree $w$ and type $\alpha=\left( \alpha_0, \ldots, \alpha_N \right)$ with an isolated singularity such that the greatest common divisor of all coefficients of $f$ is $1$.
Furthermore, let $X$ be the hypersurface in $\mathbb{P}^N_{\mathbb{Z}}\left( \alpha \right)$ defined by $f$. Let $f_p \in \mathbb{F}_p[x_0, \ldots, x_N]$ be the reduction of $f$ modulo $p$ and assume that $p \geq w(N-2)+1$.
Then $\fpt(f_p)=1=\lct(f)$ if and only if $\h\left( H^{N-1}\left( X , {\mathbb{G}}_{m,X}\right) \right)=1$.
\end{theorem}

\begin{proof}
One can show that $X$ is flat over $\mathbb{Z}$. Hence, 
Theorem 1 of \cite{Stienstra} (which also holds for quasi-homogeneous polynomials, see \cite{Yui}, section 5) yields that $H^{N-1}\left( X, {\mathbb{G}}_{m,X} \right)$ is a formal group of dimension $1$. Using the notation of Example \ref{StienstraEx} we have $J=\left\{(1, \ldots,1)\right\}$, since $d=w$. Therefore, the logarithm of the formal group law is given by 
$$l(\tau)=\mathlarger{\mathlarger{\sum}}_{m \geq 1} \frac{b_{m-1}}{m}\tau^{m},$$
where $b_{m-1}$ is the coefficient of $\left(x_0 \cdots x_n\right)^{m-1}$ in $f^{m-1}$.

Using the remark after Lemma 4.1. of \cite{Mueller} we have that $\fpt(f_p)=1$ if and only if $b_{p-1} \not \equiv 0 \mod p$. Furthermore, $b_{p-1} \not \equiv 0 \mod p$ if and only if $\ord_p \left( b_{p-1} \right) = 0$.
Finally, by Lemma \ref{VlasenkoCorThm2} it follows that $\ord_p \left( b_{p-1} \right) = 0$ if and only if $\h\left( H^{N-1}\left( X , {\mathbb{G}}_{m,X}\right) \right)=1$. 
\end{proof}

\begin{example}
Let $f=x_0^d+ \ldots + x_N^d \in \mathbb{Z}[x_0, \ldots, x_N]$ with $d=N+1$ and let $X$ be the Fermat hypersurface in $\mathbb{P}^N_{\mathbb{Z}}$ given by $f$. Let $p \geq (N+1)(N-2)+1$. Then by Theorem \ref{TheoremCYhom} we have $\fpt(f_p)=1=\lct(f)$ if and only if $\h \left( H^{N-1}\left( X , {\mathbb{G}}_{m,X}\right) \right)=1$.
\end{example}

The aim of the rest of this section is to show that a similar result as the above also holds for 
$$f=x_0^d+ \ldots + x_N^d \in \mathbb{Z}[x_0, \ldots, x_N] \text{ with } d=N+k \text{ and } k \geq 2.$$
Before we consider the case $k \geq 3$, we start with $k=2$. For this, we need the following lemma, which holds in a more general setting.

\begin{lemma}\label{Hilfslemma}
Let $f \in K[x_0, \ldots, x_N]$ be a quasi-homogeneous polynomial of degree $d$ and type $\alpha=\left( \alpha_0, \ldots, \alpha_N \right)$, where $K$ is a field of characteristic $p>0$. Let $w:=\alpha_0+ \ldots + \alpha_N$. If $wq \equiv x \mod d$ with $1 \leq x \leq w-1$ for some $q=p^e$, then $\fpt(f)< \frac{w}{d}$.
\end{lemma}

\begin{proof}
By Lemma 3.6 of \cite{Mueller} we have that $\mu_f(q) \leq \left\lceil \frac{wq-w+1}{d}\right\rceil$. By the assumption we get
$$\mu_f(q) \leq \left\lceil \frac{wq-w+1}{d}\right\rceil \leq  \frac{wq-w+1+(w-2)}{d}=\frac{wq-1}{d}.$$
Therefore, $\mu_f(q)<\frac{wq}{d}$ and $\fpt(f) \leq \frac{\mu_f(q)}{q}<\frac{w}{d}$.
\end{proof}

\begin{example}
Let $f=x_0^d+ \ldots + x_N^d \in \mathbb{Z}[x_0, \ldots, x_N]$ with $d=N+2$ and let $X$ be the Fermat hypersurface in $\mathbb{P}^N_{\mathbb{Z}}$ given by $f$.
We claim that the formal group $H^{N-1} \left( X, {\mathbb{G}}_{m,X} \right)$ is the direct sum of $N+1$ copies of a $1$-dimensional formal group law $F$ and that $\fpt(f_p)=\lct(f)=\frac{N+1}{N+2}$ if and only if $\h(F)=1$.

For the proof of this, we use the notation of Example \ref{StienstraEx} and compute $$J=\left\{(2,1,\ldots,1), (1,2,1,\ldots,1), \ldots, (1, \ldots, 1,2)\right\}.$$
For $i,j \in J$ we denote by $b_{m-1,i,j}$ the coefficient of $x_0^{mj_0-i_0} \ldots x_N^{mj_N-i_N}$ in $f^{m-1}$. It is an easy computation to see that $b_{m-1,i,j}=0$ if $i \neq j$.
Therefore, the formal group $H^{N-1} \left( X, {\mathbb{G}}_{m,X} \right)$ is the direct sum of $N+1$ copies of the $1$-dimensional formal group law $F(\tau,\eta)=l^{-1}\left( l(\tau)+l(\eta) \right)$ with logarithm
$$l\left(\tau\right)=\mathlarger{\mathlarger{\sum}}_{n \geq 0} \frac{1}{nd+1} \frac{(nd)!}{n!^{N}(2n)!} \tau^{nd+1}
=\mathlarger{\mathlarger{\sum}}_{m \geq 1} \frac{b_{m-1}}{m} \tau^{m},$$
where
$$b_{m-1}= \begin{cases} \frac{(nd)!}{n!^{N}(2n)!}, &\mbox{if  $m=1+nd$ for some $n \in \mathbb{Z}$} \\
0, & \mbox{otherwise.} \end{cases}$$

Now, we show that $\fpt(f_p)=\lct(f)=\frac{N+1}{N+2}$ if and only if $\h(F)=1$.

For this, remark that by Lemma \ref{VlasenkoCorThm2} it follows that $\h(F)=1$ if and only if $\ord_p \left( b_{p-1} \right) = 0$, which is equivalent to $b_{p-1} \not \equiv 0 \mod p$.
For $p \equiv 1 \mod d$, we have that $nd=p-1$ and $2n=\frac{2(p-1)}{d}$ are smaller than $p$.
Since $\ord_p(s!)=\sum_{i \geq 1} \left\lfloor \frac{s}{p^i} \right\rfloor$ for all $s \in \mathbb{N}$, this means that the $p$-adic valuation of $\frac{(nd)!}{n!^{N}(2n)!}$ is zero. Therefore, $b_{p-1} \not \equiv 0 \mod p$ if and only if $p \equiv 1 \mod d$.

Hence, it remains to prove that $p \equiv 1 \mod d$ is equivalent to $\fpt(f_p)=\lct(f)$. Using Example 4.2. of \cite{MTW} one gets that $\fpt(f_p)=\lct(f)$ if $p \equiv 1 \mod d$. Now let $\fpt(f_p)=\lct(f)$. By Lemma \ref{Hilfslemma} we conclude that $(N+1)q \equiv x \mod N+2$ with $x \in \left\{0,N+1\right\}$ for all $q$. If $x=0$ for $q=p$ we get $(N+1)p \equiv 0 \mod N+2$. Therefore $p \equiv 0 \mod N+2$, which is a contradiction. If $x=N+1$ one gets $(N+1)p \equiv N+1 \mod N+2$, hence $p \equiv 1 \mod N+2$.
\end{example}

Next, we consider the general case $d=N+k$ with $k \geq 3$.

\begin{lemma}\label{heightFermat}
Let $f=x_0^d+ \ldots + x_N^d \in \mathbb{Z}[x_0, \ldots, x_N]$ with $d=N+k$ for $k \geq 3$ and such that $N \geq 2(k-1)$ and let $X$ be the hypersurface in $\mathbb{P}^N_{\mathbb{Z}}$ given by $f$. 
Then the formal group $H^{N-1} \left( X, {\mathbb{G}}_{m,X} \right)$ is the direct sum of $\binom{N+k-1}{N}$ $1$-dimensional formal groups, which are all of height $1$ if and only if $p \equiv 1 \mod d$.
\end{lemma}

\begin{proof}
As in Example \ref{StienstraEx} let 
$$J=\left\{i=\left( i_0, \ldots, i_N\right) \in \mathbb{Z}^{N+1}| i_0, \ldots, i_N \geq 1 \text{ and } i_0+\ldots+i_N=d\right\}$$
and for $i, j \in J$ let $b_{m-1,i,j}$ be the coefficient of $x_0^{mj_0-i_0} \cdots x_N^{mj_N-i_N}$ in $f^{m-1}$. 
We prove the lemma via the following steps:

(1) We show, that for $i \neq j$ one has $b_{m-1,i,j}=0$:
Since $k \geq 2$ and $N \geq 2(k-1)$ it follows that $N \geq k$.
The elements of the set $J$ are tuples $i=(i_0, \ldots, i_N)$ with $i_0, \ldots, i_N \geq 1$ and 
\begin{align*} \label{allgFermatAbsch}
i_0+ \ldots+ i_N=d=N+k=(N+1)+(k-1) \leq (N+1)+(N-1),
\end{align*}
i.e. each entry $i_n$ is at least one and further $k-1\leq N-1$ has to be distributed in the entries of $i$.

Since $N \geq 2(k-1)$, it follows that $\frac{N+1}{2}>k-1$, 
 i.e. more than half of the entries of a tuple $i \in J$ are equal to $1$. This means that if $j=(j_0, \ldots, j_N) \in J$ is a second tuple, then there exists at least one position $s$ with $i_s=1=j_s$.

Now write
$$f^{m-1}=\mathlarger{\mathlarger{\sum}}_{\beta_0+\ldots+\beta_N=m-1} \binom{m-1}{\beta_0, \ldots, \beta_N} x_0^{d\beta_0} \cdots x_N^{d \beta_N}.$$
Then we have
\begin{align*}
d \beta_n &= mj_n-i_n \text{ for all } n \neq s\\
d \beta_s &= m-1.
\end{align*}
The last equality shows that $m \equiv 1 \mod d$ and the first equality then yields
$0 \equiv mj_n-i_n \mod d \equiv j_n-i_n \mod d$, i.e. $j_n \equiv i_n \mod d$ for all $n \neq s$. But since $i_n, j_n \leq N$ and $d=N+k > N$ it follows that $i_n=j_n$ for all $n \neq s$ and therefore $i=j$.

(2) Part (1) of this proof means, that the logarithm $l \left( \tau \right)$ of the formal group $H^{N-1} \left( X, {\mathbb{G}}_{m,X} \right)$ of dimension $\#J=\binom{N+k-1}{N}$ is given by $\left(l_i \left( \tau_i \right)\right)_{i \in J}$, where
$$l_i\left(\tau_i\right)=\mathlarger{\mathlarger{\sum}}_{m \geq 1} \frac{b_{m-1,i,i}}{m} \tau_i^{m}.$$
and one can compute that
$$b_{m-1,i,i}= \begin{cases} \frac{(ad)!}{(ai_0)!\cdots (ai_N)!}, &\mbox{if  $m=1+ad$ for some $a \in \mathbb{Z}$} \\
0, & \mbox{otherwise.} \end{cases}$$

(3) By (1) and (2) we know that $H^{N-1} \left( X, {\mathbb{G}}_{m,X} \right)$ is the direct sum of $\binom{N+k-1}{N}$ formal groups $\left( F_i \right)_{i \in J}$, where $F_i(\tau_i, \eta_i)=l_i^{-1} \left( l_i(\tau_i) + l_i(\eta_i)\right)$.
We prove that $\h(F_i)=1$ for all $i \in J$ if and only if $p \equiv 1 \mod d$.
For this, Lemma \ref{VlasenkoCorThm2} shows that $\h(F_i)=1$ if and only if $b_{p-1,i,i} \not \equiv 0 \mod p$. For $p \equiv 1 \mod d$, we have that $ad=p-1$ and $a i_r=\frac{p-1}{d} i_r < p-1$ for all $0 \leq r \leq N$ and hence the $p$-adic valuation of $\frac{(ad)!}{(ai_0)!\cdots (ai_N)!}$ is zero. Therefore, it follows that $b_{p-1,i,i} \not\equiv  0 \mod p$ if and only if $p \equiv 1 \mod d$.
\end{proof}

The following lemma computes the $F$-pure threshold of Fermat hypersurfaces.

\begin{lemma}\label{FthresholdFermat}
Let $f=x_0^d+ \ldots + x_N^d \in \mathbb{Z}[x_0, \ldots, x_N]$ with $d=N+k$ for $k \geq 2$ and such that $N > k-2$. Furthermore, let $d \not \equiv 0 \mod p$. Then $p \equiv 1 \mod d$ if and only if $\fpt(f_p)=\lct(f)=\frac{N+1}{d}$.
\end{lemma}

\begin{proof}
First, let $p \equiv 1 \mod d$. Then by example 4.2 of \cite{MTW} it follows
that $\fpt(f_p)=\lct(f)$.

Now, we show that if $p \not \equiv 1 \mod d$, then $\fpt\left(f_p \right) < \lct(f)$. For this, remember that $\fpt(f_p)= \lim\limits_{e \rightarrow \infty} \frac{\mu_{f_p}(p^e)}{p^e}$, where $\mu_{f_p}(p^e)=\min \left\{n \in \mathbb{N} \big| f_p^n \in \mathfrak{m}^{[p^e]}\right\}$ and 
$$f_p^n=\mathlarger{\mathlarger{\sum}}_{\beta_0+ \ldots + \beta_N=n}  \binom{n}{\beta_0, \ldots, \beta_N} x_0^{d \beta_0} \cdots x_N^{d \beta_N}.$$
We claim that $\frac{\mu_{f_p}(p)}{p} < \frac{N+1}{d}$. Once we have shown this, it follows that $\fpt(f_p) \leq \frac{\mu_{f_p}(p)}{p} < \frac{N+1}{d}$.
In order to show $\frac{\mu_{f_p}(p)}{p} < \frac{N+1}{d}$ or equivalently $\mu_{f_p}(p)< \frac{p(N+1)}{d}$, it is enough to show that $f_p^{\left\lfloor \frac{p(N+1)}{d}\right\rfloor} \in \mathfrak{m}^{[p]}$, since $d \not \equiv 0 \mod p$.
For this, it is enough to show that 
\begin{equation*}\label{round}
d \left\lceil \frac{\left\lfloor \frac{p(N+1)}{d}\right\rfloor}{N+1}\right\rceil \geq p,
\end{equation*}
We now consider the following two cases:\\
\textbf{Case 1}: $N+1 \nmid \left\lfloor \frac{p(N+1)}{d}\right\rfloor$\\
Clearly one has $\left\lfloor \frac{p(N+1)}{d}\right\rfloor+1 \geq \frac{p(N+1)}{d}$, hence $\frac{\left\lfloor \frac{p(N+1)}{d}\right\rfloor}{N+1}+\frac{1}{N+1} \geq \frac{p}{d}$. Since $N+1$ does not divide $\left\lfloor \frac{p(N+1)}{d}\right\rfloor$ by assumption, this last inequality yields 
$\left\lceil \frac{\left\lfloor \frac{p(N+1)}{d}\right\rfloor}{N+1}\right\rceil \geq \frac{p}{d}$.\\
\textbf{Case 2}: $N+1 \mid \left\lfloor \frac{p(N+1)}{d}\right\rfloor$\\
Write $p=\lambda d +r$ with $0 \leq r<d$ and $r \neq 1$, since $p \not \equiv 1 \mod d$. Thus $N+1$ divides
$$\left\lfloor \frac{p(N+1)}{d}\right\rfloor= \left\lfloor \frac{(\lambda d+r) (N+1)}{d}\right\rfloor=\left\lfloor \lambda (N+1)+\frac{r(N+1)}{d}\right\rfloor=  \lambda (N+1)+\left\lfloor\frac{r(N+1)}{d}\right\rfloor.$$
Since $\frac{r}{d}<1$, it follows that $\frac{r(N+1)}{d}<N+1$ and since $\left\lfloor\frac{r(N+1)}{d}\right\rfloor$ must be divisible by $N+1$ we conclude that $\left\lfloor\frac{r(N+1)}{d}\right\rfloor=0$.
This means that $\frac{r(N+1)}{d}<1$ or equivalently $r < \frac{d}{N+1}=\frac{N+k}{N+1}$. Since $k < N+2$, we have
$r < \frac{N+k}{N+1}< \frac{2N+2}{N+1}=2$. By assumption $r \neq 1$, hence $r=0$ and $p=d$, since $p$ is a prime.
But this is a contradiction to our assumptions.
\end{proof}

Combining Lemma \ref{heightFermat} and Lemma \ref{FthresholdFermat} we obtain:

\begin{corollary}\label{CorFermat}
Let $f=x_0^d+ \ldots + x_N^d \in \mathbb{Z}[x_0, \ldots, x_N]$ with $d=N+k$ for $k \geq 3$ and such that $N \geq 2(k-1)$ and let $X$ be the hypersurface in $\mathbb{P}^N_{\mathbb{Z}}$ given by $f$. 
Furthermore, let $d \not \equiv 0 \mod p$.
Then $H^{N-1} \left( X, {\mathbb{G}}_{m,X} \right)$ is a direct sum of formal groups  
of dimension $1$, which are all of height $1$ if and only if $\fpt(f_p)=\lct(f)=\frac{N+1}{d}$.
\end{corollary}

If one combines this with the result of Koblitz in \cite{Koblitz}, one obtains that the two conditions above are also equivalent to the Frobenius action on $H^{N-1} \left(X, \mathcal{O}_X \right)$ being bijective.

In the proofs of Corollary \ref{CorFermat} and Theorem \ref{TheoremCYhom} we have seen that the $F$-pure threshold is equal to the log canonical threshold if and only if the height of the corresponding Artin-Mazur formal group is equal to its dimension. Or, equivalently, since $\fpt(f_p) \leq \lct(f)$ for all $p$ it means that the $F$-pure threshold is equal to its greatest possible value if and only if the height is equal to its smallest possible value (see Lemma \ref{direkteSummefG}).
Since we did not find any counterexample for this so far, this leads us to suspect that this could be the case for all quasi-homogeneous polynomials.

\section{Counterexamples}

Let $R:=K[x_0, \ldots, x_N]$ be the graded polynomial ring with $\alpha_i:=\deg(x_i)$ over an algebraically closed field $K$ of characteristic $p>0$.
Let $f \in R$ be a quasi-homogeneous polynomial of degree $w= \alpha_0+ \ldots + \alpha_N$ and type $\alpha=\left( \alpha_0, \ldots, \alpha_N\right)$ with an isolated singularity.
Theorem 3.9 together with Theorem 5.1 of \cite{Mueller} yield that 
$$\fpt(f)=1-\frac{a}{p}$$ 
with $0 \leq a \leq N-1$ for $p \geq w(N-2)+1$, where $a$ is the order of vanishing of the Hasse invariant on a certain deformation space of $X=\Proj(R/fR) \subset \mathbb{P}^N \left(\alpha \right)$.
Theorem 3.1 shows that $a=0$ if and only if $\h \left( H^{N-1}\left( X , {\mathbb{G}}_{m,X}\right) \right)=1$.

Therefore, one may ask whether the other possible values of the $F$-pure threshold (i.e. $1 \leq a \leq N-1$) can also be characterized by $\h\left(H^{N-1}\left( X , {\mathbb{G}}_{m,X}\right)\right)$.
However, in this section we will give two examples of weighted Delsarte $K3$ surfaces which show that the answer is negative. 

First, let us briefly recall the definition of a weighted Delsarte $K3$ surface. For more details, we refer the reader to \cite{Goto}.
Let $N=3$ and assume that 
\begin{align*}
p \nmid \alpha_i \text{ for } 0 \leq i \leq 3 \text{ and }
\end{align*}
\begin{align*}
\ggT(\alpha_i, \alpha_j, \alpha_k)=1 \text{ for all } \left\{i, j, k\right\} \subset \left\{0, 1, 2, 3\right\}.
\end{align*}
Let $m$ be a positive integer such that $p \nmid m$ and let $A=(a_{ij}) \in \mathbb{Z}^{4 \times 4}$ be a matrix 
such that
\begin{enumerate}
	\item $a_{ij} \geq 0$ and $p \nmid a_{ij}$ for all $(i, j)$,
	\item given $j$ there is some $i$, such that $a_{ij}=0$,
	\item $p \nmid \det(A)$,
	\item $\sum\limits_{j=0}^3 \alpha_j a_{ij}=m$ for all $0 \leq i \leq 3$, i.e. $A \alpha^T = \begin{pmatrix} m, m, m, m \end{pmatrix}^T$.
\end{enumerate}
 
A \textbf{weighted Delsarte surface} in $\mathbb{P}^3(\alpha)$ of degree $m$ with matrix $A$ is defined to be the surface $X_A \subset \mathbb{P}^3(\alpha)$ given by
$$\sum_{i=0}^3 x_0^{a_{i0}} x_1^{a_{i1}} x_2^{a_{i2}} x_3^{a_{i3}}=0.$$

Let $p: \mathbb{A}^{4} \setminus \left\{0\right\} \rightarrow \mathbb{P}^3(\alpha)$ be the canonical projection. Then the scheme closure of $p^{-1}(X_A)$ in $\mathbb{A}^{4}$ is called the affine quasicone over $X_A$.
We say that $X_A \subset \mathbb{P}^3(\alpha)$ is \textbf{quasi-smooth}, if its affine quasicone is smooth outside the origin (see \cite{Dolgachev}).
Furthermore, we say that $X_A$ is \textbf{in general position relative to $\mathbb{P}^3(\alpha)_{\text{sing}}$} if $\codim_{X}(X \cap \mathbb{P}^3(\alpha)_{\text{sing}}) \geq 2$, where $\mathbb{P}^3(\alpha)_{\text{sing}}$ denotes the singular locus of $\mathbb{P}^3(\alpha)$ (see \cite{Gotosymplectic}).

Weighted Delsarte surfaces are in general singular surfaces.
If $X_A$ is quasi-smooth and in general position relative to $\mathbb{P}^3(\alpha)_{\text{sing}}$, then the minimal resolution $\widetilde{X_A}$ of $X_A$ is a $K3$ surface if and only if $m=\alpha_0+\alpha_1+\alpha_2+\alpha_3$.
If this is the case, then we call $X_A$ a \textbf{weighted Delsarte $K3$ surface} in $\mathbb{P}^3(\alpha)$ of degree $m$ with matrix $A$.

Let 
$$e_A:=\frac{\left|\det(A)\right|}{g},$$
where $g$ is the $\gcd$ of all column sums of the adjugate matrix of $A$ and of $\left|\det(A)\right|$. Goto gives the following criterion for the formal Brauer group of $\widetilde{X_A}$ to have infinite height. 

\begin{lemma}[{\cite[Proposition 2.2 \& Remark 2.1]{Goto}}]\label{supersingular}
Let $X_A$ be a weighted Delsarte $K3$ surface with matrix $A$. Then the height of the formal Brauer group of the minimal resolution $\widetilde{X_A}$ of $X_A$ is infinite (i.e. $\widetilde{X_A}$ is supersingular) if and only if $p^\mu \equiv -1 \mod e_A$ for some integer $\mu \geq 1$. 
\end{lemma}

Furthermore, he explains how to compute the height of the formal Brauer group of a weighted Delsarte $K3$ surface if it is finite:

\begin{theorem}[{{\cite[Theorem 3.2]{Goto}}}]\label{finiteheight}
Let $X_A$ be a weighted Delsarte $K3$ surface with  matrix $A$. Assume that there is no integer $\mu \geq 1$ such that $p^\mu \equiv -1 \mod e_A$. Then the height of the formal Brauer group of the minimal resolution $\widetilde{X_A}$ of $X_A$ is equal to the order of $p$ modulo $e_A$.
\end{theorem}

We use these two results to give two examples of weighted Delsarte $K3$ surfaces. The first one, will have the same height but different $F$-pure threshold for varying $p$ and the second one will have the same $F$-pure threshold but the height will differ for two different primes $p$.

\begin{example}
Assume that $p \neq 2,5$. Consider $f=x^2+y^5+z^5+w^{10} \in K[x,y,z,w]$, which is quasi-homogeneous of degree $10$ and weight $(5,2,2,1)$.
Let $X_A$ be the weighted Delsarte surface in $\mathbb{P}^3(5,2,2,1)$ defined by $f$, i.e. defined by the matrix
$$A=\begin{pmatrix} 
2 & 0 & 0 & 0\\ 
0 & 5 & 0 & 0\\ 
0 & 0 & 5 & 0\\ 
0 & 0 & 0 & 10
\end{pmatrix}.$$
Using the methods of \cite{Goto} we computed the height of the formal Brauer group of the minimal resolution $\widetilde{X_A}$ of $X_A$.
Since $\sqrt{J(f)}=(x,y,z,w)$ and $\codim_{X_A}(X_A \cap \mathbb{P}^3(\alpha)_{sing}) \geq 2$, $X_A$ is quasi-smooth and in general position relative to $\mathbb{P}^3(\alpha)_{sing}$. Furthermore, $m=10=\alpha_0 + \alpha_1 + \alpha_2 + \alpha_3$ and therefore the minimal resolution $\widetilde{X_A}$ of $X_A$ is a $K3$ surface.
One has $\det(A)=500$, $g=50$ and therefore $e_A=10$. 
Thus, Lemma \ref{supersingular}
shows that the height of the formal Brauer group of $\widetilde{X_A}$ is infinite if and only if there exists some $\mu \geq 1$ such that $p^{\mu} \equiv -1 \mod 10$, i.e. $p \equiv 3, 7, 9 \mod 10$.  

Using the PosChar-package of Macaulay 2 \cite{PosChar} we also computed the $F$-pure threshold of $f$. We obtained the following results:
\renewcommand{\arraystretch}{1.3}
\begin{center}
\begin{tabular}{|c|c|c|c|c|c|c|} \hline
$p$ & $3$ & $7$ & $11$ & $13$ & $17$ & $19$ \\ \hline
$\h$ & $\infty$ & $\infty$ & $1$ & $\infty$ & $\infty$ & $\infty$ \\ \hline
$\fpt(f)$ & $1-\frac{1}{p}$ & $1-\frac{1}{p}$ & $1$ & $1-\frac{1}{p}$ & $1-\frac{1}{p}$ & $1-\frac{2}{p}$ \\ \hline
 \end{tabular}
\end{center}
In particular, one can see that for $p=17$ and $p=19$ the height is the same but the $F$-pure threshold is different.
\end{example}

\begin{example}
Assume that $p \neq 2,3$. Consider $f=x^8y+y^6z+z^3+xw^2 \in K[x,y,z,w]$, which is quasi-homogeneous of degree $9$ and weight $(1,1,3,4)$.
Let $X_A$ be the weighted Delsarte surface in $\mathbb{P}^3(1,1,3,4)$ defined by $f$, i.e. defined by the matrix
$$A=\begin{pmatrix} 
8 & 1 & 0 & 0\\ 
0 & 6 & 1 & 0\\ 
0 & 0 & 3 & 0\\ 
1 & 0 & 0 & 2
\end{pmatrix}.$$
To compute the height of the formal Brauer group of the minimal resolution $\widetilde{X_A}$ of $X_A$ one checks that $\sqrt{J(f)}=(x,y,z,w)$ and $\codim_{X_A}(X_A \cap \mathbb{P}^3(\alpha)_{sing}) \geq 2$, so $X_A$ is quasi-smooth and in general position relative to $\mathbb{P}^3(\alpha)_{sing}$. Furthermore, $m=9=\alpha_0 + \alpha_1 + \alpha_2 + \alpha_3$ and therefore the minimal resolution $\widetilde{X_A}$ of $X_A$ is a $K3$ surface.
We compute that $\det(A)=288$, $g=9$ and therefore $e_A=32$.
Using Theorem \ref{finiteheight} we get that the height of the formal Brauer group of $\widetilde{X_A}$ is given by
$$\h= \begin{cases} 1, &\mbox{if  $p \equiv 1 \mod 32$} \\
2, &\mbox{if  $p \equiv \pm 15 \mod 32$} \\
4, &\mbox{if  $p \equiv \pm 7, \pm 9 \mod 32$} \\
8, &\mbox{if  $p \equiv \pm 3, \pm 5, \pm 11, \pm 13 \mod 32$}.
 \end{cases}$$

Combined with the $F$-pure thresholds of $f$ one obtains:
\renewcommand{\arraystretch}{1.3}
\begin{center}
\begin{tabular}{|c|c|c|c|c|c|c|c|} \hline
$p$ & $3$ & $5$ & $7$ & $11$ & $13$ & $17$ & $19$ \\ \hline
$\h$ & $8$ & $8$ & $4$ & $8$ & $8$ & $2$ & $8$ \\ \hline
$\fpt(f)$ & $1-\frac{1}{p}$ & $1-\frac{1}{p}$ & $1-\frac{1}{p}$ & $1-\frac{1}{p}$ & $1-\frac{1}{p}$ & $1-\frac{1}{p}$ & $1-\frac{1}{p}$  \\ \hline
 \end{tabular}
\end{center}
In particular, in this case the $F$-pure threshold is $1-\frac{1}{p}$ for all $p$, but the height differs.
\end{example}


\newcommand{\etalchar}[1]{$^{#1}$}

\end{document}